# A RESHETNYAK-TYPE LOWER SEMICONTINUITY RESULT FOR LINEARISED ELASTO-PLASTICITY COUPLED WITH DAMAGE IN $W^{1,n}$

VITO CRISMALE AND GIANLUCA ORLANDO

ABSTRACT. In this paper we prove a lower semicontinuity result of Reshetnyak type for a class of functionals which appear in models for small-strain elasto-plasticity coupled with damage. To do so we characterise the limit of measures $\alpha_k \operatorname{E} u_k$ with respect to the weak convergence $\alpha_k \rightharpoonup \alpha$ in $W^{1,n}(\Omega)$ and the weak${}^*$ convergence $u_k \stackrel{*}{\rightharpoonup} u$ in $BD(\Omega)$, E denoting the symmetrised gradient. A concentration compactness argument shows that the limit has the form $\alpha \operatorname{E} u + \eta$, with $\eta$ supported on an at most countable set.

## Contents

## 1. Introduction

In this paper we prove a lower semicontinuity result of Reshetnyak type for a class of functionals which appear in models for small-strain elasto-plasticity coupled with damage. The functionals $\mathcal{H}(\alpha, p)$ that we consider depend on Sobolev functions $\alpha$, the damage variables, and on bounded Radon measures $p$, the plastic strains.

In small-strain plasticity, the linearized strain $\operatorname{E} u$, defined as the symmetric part of the spatial gradient of the displacement $u \colon \Omega \to \mathbb{R}^n$, is decomposed as the sum of the elastic strain $e \in L^2(\Omega; \mathbb{M}^{n \times n}_{\mathrm{sym}})$, and of the plastic strain $p \in \mathcal{M}_b(\Omega; \mathbb{M}^{n \times n}_{\mathrm{sym}})$, i.e., $p$ is a bounded Radon measure with values in the space of symmetric matrices $\mathbb{M}^{n \times n}_{\mathrm{sym}}$. In perfect plasticity (without damage), the energy dissipated in the evolution of the plastic strain is described in terms of the so-called plastic potential, defined in accordance to the theory of convex functions of measures by

$$\int_\Omega H\Big(\frac{\mathrm{d} p}{\mathrm{d} |p|}(x)\Big) \, \mathrm{d}|p|(x), \quad \text{for } p \in \mathcal{M}_b(\Omega; \mathbb{M}^{n \times n}_{\mathrm{sym}}).$$

In the formula above, $\mathrm{d} p / \mathrm{d} |p|$ is the Radon-Nikodym derivative of $p$ with respect to its total variation $|p|$ and $H$ is the support function of a set $K + \mathbb{R} I$, $I$ being the identity matrix and $K$ the convex compact set of the space of $n \times n$ trace-free matrices where the deviatoric part of the stress is constrained to lie. In particular, $H \colon \mathbb{M}^{n \times n}_{\mathrm{sym}} \to [0, \infty]$ is convex, lower semicontinuous, and positively 1-homogeneous. We refer to [12] for all the details about the mathematical formulation of small-strain perfect plasticity.

In presence of damage, the constraint set depends on the real-valued damage variable $\alpha$. Here we assume a multiplicative dependence, that is $K(\alpha) = V(\alpha)K$, with $V\colon \mathbb{R} \to [0,\infty)$ lower semicontinuous. In this setting the plastic potential becomes

$$\mathcal{H}(\alpha, p) := \int_\Omega V(\alpha(x)) H\Big(\frac{\mathrm{d}p}{\mathrm{d}|p|}(x)\Big)\, \mathrm{d}|p|(x)\,. \tag{1.1}$$

The functional above is sequentially lower semicontinuous with respect to the uniform convergence in $\alpha$ and the weak* convergence in $p$, as a consequence of Reshetnyak's Lower Semicontinuity Theorem (see, e.g., [6, Theorem 2.38]).

The lower semicontinuity of the plastic potential is, in general, a major difficulty in small-strain plasticity when the constraint set depends on an additional variable. For instance, in non-associative plasticity (cf. [13, 7] and the recent [18]) such variable lacks continuity, and Reshetnyak's Theorem cannot be applied directly. The way out consists in replacing the original additional variable by a mollified one.

In gradient damage models, the total energy features a term in $\nabla \alpha$ which provides uniform bounds for $\alpha$ in $W^{1,q}(\Omega)$, for a suitable $q > 1$. When one considers the coupling with plasticity in the case $q > n$, the functional in (1.1) is defined by choosing the continuous representative of $\alpha$ and is sequentially lower semicontinuous with respect to the weak convergence in $W^{1,q}(\Omega)$, in view of the compact embedding of $W^{1,q}(\Omega)$ in $\mathcal{C}(\overline{\Omega})$. In particular, the minimum problems involved in the variational approach to the existence of quasistatic evolutions admit solutions, cf. [8].

However, in many mechanical models [27, 28, 20, 3, 24, 4, 23] the natural space for the damage variable is $W^{1,q}(\Omega)$ for some exponent $q \leq n$, usually the Hilbert space $H^1(\Omega)$. Here we focus our attention on the critical case $q = n$, which in particular covers two dimensional models with damage in $H^1(\Omega)$. Observe that a function $\alpha \in W^{1,n}(\Omega)$ does not always admit a continuous representative. Nonetheless, the precise representative $\widetilde{\alpha}$ of $\alpha$ is defined up to a set of $n$-capacity zero. In particular, this exceptional set has $\mathcal{H}^{n-1}$-measure zero and thus it is $|p|$-negligible. The functional in (1.1) is therefore well-defined upon choosing this precise representative $\widetilde{\alpha}$.

The main result in this paper is the following.

**Theorem 1.1.** *Assume that $\Omega$ is a bounded, open set with Lipschitz boundary. Let $V\colon \mathbb{R} \to [0,\infty]$ and $H\colon \Omega \times \mathbb{M}^{n \times n}_{\mathrm{sym}} \to [0,\infty]$ be lower semicontinuous functions. Assume that $H$ is positively 1-homogeneous and convex in the second variable. Let $\alpha_k, \alpha \in W^{1,n}(\Omega)$, $u_k, u \in BD(\Omega)$, $e_k, e \in L^q(\Omega; \mathbb{M}^{n \times n}_{\mathrm{sym}})$ for some $q > 1$, and $p_k, p \in \mathcal{M}_b(\Omega; \mathbb{M}^{n \times n}_{\mathrm{sym}})$. Assume that*

$$\mathrm{E}u_k = e_k + p_k \quad \textit{in } \mathcal{M}_b(\Omega; \mathbb{M}^{n \times n}_{\mathrm{sym}})\,, \tag{1.2}$$

$$u_k \overset{*}{\rightharpoonup} u \quad \textit{weakly* in } BD(\Omega)\,, \tag{1.3}$$

$$e_k \rightharpoonup e \quad \textit{weakly in } L^q(\Omega; \mathbb{M}^{n \times n}_{\mathrm{sym}})\,, \tag{1.4}$$

$$\alpha_k \rightharpoonup \alpha \quad \textit{weakly in } W^{1,n}(\Omega)\,.$$

*Then*

$$\int_\Omega V(\widetilde{\alpha}(x))\, H\Big(x, \frac{\mathrm{d}p}{\mathrm{d}|p|}(x)\Big)\, \mathrm{d}|p|(x) \leq \liminf_{k \to +\infty} \int_\Omega V(\widetilde{\alpha}_k(x))\, H\Big(x, \frac{\mathrm{d}p_k}{\mathrm{d}|p_k|}(x)\Big)\, \mathrm{d}|p_k|(x)\,. \tag{1.5}$$

To illustrate the proof of Theorem 1.1, we consider now the simplified case $V(\alpha) = \alpha$, $0 \leq \alpha \leq 1$, and $H(x, \xi) = |\xi|$. The starting point is the following Leibniz formula (Proposition 3.5)

$$\widetilde{\alpha}_k\, \mathrm{E}u_k = \mathrm{E}(\widetilde{\alpha}_k\, u_k) - \nabla \alpha_k \odot u_k\,,$$

where $\odot$ denotes the symmetric tensor product. If the sequence $u_k$ were bounded in $L^\infty(\Omega; \mathbb{R}^n)$, then $\nabla \alpha_k \odot u_k$ would converge weakly in $L^n(\Omega; \mathbb{M}^{n \times n}_{\mathrm{sym}})$ to $\nabla \alpha \odot u$, and the formula above would easily imply that $\widetilde{\alpha}_k\, \mathrm{E}u_k \overset{*}{\rightharpoonup} \widetilde{\alpha}\, \mathrm{E}u$. In different contexts where truncation arguments are allowed, this makes possible to prove the lower semicontinuity of the plastic potential, cf. [14, Proposition 2.3] for the antiplane setting of this model and [9, Theorem 3.1] for the coupling of damage and strain gradient plasticity.

Here we are able to give a precise description of the weak* limit of the sequence $\widetilde{\alpha}_k\, \mathrm{E}u_k$, which may differ from $\widetilde{\alpha}\, \mathrm{E}u$ (cf. Example 3.1). Specifically, a concentration compactness argument in the spirit of [21] allows us to prove in Theorem 3.2 that $\widetilde{\alpha}_k\, \mathrm{E}u_k \overset{*}{\rightharpoonup} \widetilde{\alpha}\, \mathrm{E}u + \eta$, where $\eta$ is a measure concentrated on an at most countable set. In particular, $\widetilde{\alpha}_k\, p_k \overset{*}{\rightharpoonup} \widetilde{\alpha}\, p + \eta$. Passing to the total variations, this entails the desired lower semicontinuity since $\widetilde{\alpha}\, p$ and $\eta$ are mutually singular.

We stress that this type of proof only works in the critical case $\alpha \in W^{1,q}(\Omega)$ with $q = n$. Indeed, Example 3.7 shows that if $q < n$, it may happen that $\widetilde{\alpha}_k\, \mathrm{E}u_k \overset{*}{\rightharpoonup} \widetilde{\alpha}\, \mathrm{E}u + \eta$, where $\eta$ is not singular with respect to $\widetilde{\alpha}\, \mathrm{E}u$. The case $q < n$ will be the subject of a future investigation. We remark that when $H(\xi) = |\xi|$ and $e_k \to e$ strongly in $L^2(\Omega; \mathbb{M}^{n \times n}_{\mathrm{sym}})$, the plastic potential is lower semicontinuous even in the case $q < n$, as proven in [10, Section 4.6]. Indeed, these conditions on $H$ and $e_k$ allow for a slicing argument as in [15] which reduces the proof to the one-dimensional setting. This technique is however not suited to the case where $e_k$ is only a weakly convergent sequence.

The paper is structured as follows. In Section 2 we fix the notation and we collect some preliminary results. Section 3 is devoted to the study of the weak* limit of sequences $\widetilde{\alpha}_k\, \mathrm{E}u_k$: there we provide some explicit examples of concentration effects and we prove that the excess measure in the limit is concentrated on an at most countable set. Section 4 contains the proof of Theorem 1.1. Finally, in Section 5 we apply Theorem 1.1 to show the esistence of energetic solutions for a model which couples small-strain plasticity and damage in $W^{1,n}(\Omega)$.

## 2. Notation and preliminary results

**Notation.** Throughout the paper we assume that $n \geq 2$. The Lebesgue measure in $\mathbb{R}^n$ is denoted by $\mathcal{L}^n$, while $\mathcal{H}^s$ is the $s$-dimensional Hausdorff measure.

The space of $n \times n$ symmetric matrices is denoted by $\mathbb{M}^{n \times n}_{\mathrm{sym}}$; it is endowed with the euclidean scalar product $A : B := \mathrm{tr}(AB^T)$, and the corresponding euclidean norm $|A| := (A : A)^{1/2}$. The symmetrised tensor product $a \odot b$ of two vectors $a, b \in \mathbb{R}^n$ is the symmetric matrix with components $(a_i b_j + a_j b_i)/2$.

**Measures.** Let $\Omega$ be an open set in $\mathbb{R}^n$. The space of bounded $\mathbb{R}^m$-valued Radon measures is denoted by $\mathcal{M}_b(\Omega; \mathbb{R}^m)$. This space can be regarded as the dual of the space $\mathcal{C}_0(\Omega; \mathbb{R}^m)$ of $\mathbb{R}^m$-valued continuous functions on $\overline{\Omega}$ vanishing on $\partial \Omega$. The notion of weak* convergence in $\mathcal{M}_b(\Omega; \mathbb{R}^m)$ refers to this duality. Moreover, we denote by $\mathcal{M}_b^+(\Omega)$ the space of non-negative bounded Radon measures. If $f \in L^1(\Omega; \mathbb{R}^m)$, we shall always identify the bounded Radon measure $f\mathcal{L}^n$ with the function $f$.

Let us consider a lower semicontinuous function $H \colon \Omega \times \mathbb{R}^m \to [0, \infty]$, positively 1-homogeneous and convex in the second variable and let us consider the functional defined in accordance to the theory of convex functions of measures

$$\mathrm{H}(\mu) := \int_\Omega H\Big(x, \frac{\mathrm{d}\mu}{\mathrm{d}|\mu|}(x)\Big)\, \mathrm{d}|\mu|(x), \quad \text{for } \mu \in \mathcal{M}_b(\Omega; \mathbb{R}^m),$$

where $\mathrm{d}\mu/\mathrm{d}|\mu|$ is the Radon-Nikodym derivative of $\mu$ with respect to its total variation $|\mu|$.

*Remark* 2.1. Let $\mu, \nu \in \mathcal{M}_b(\Omega; \mathbb{R}^m)$. If $|\mu|$ and $|\nu|$ are mutually singular, then $\mathrm{H}(\mu + \nu) = \mathrm{H}(\mu) + \mathrm{H}(\nu)$ (cf. [6, Proposition 2.37]).

We recall the classical Reshetnyak's Lower Semicontinuity Theorem [29]. For a proof we refer to [6, Theorem 2.38].

**Theorem 2.2** (Reshetnyak's Lower Semicontinuity Theorem)**.** *Let $\Omega$ be an open set in $\mathbb{R}^n$. Let $\mu_k, \mu \in \mathcal{M}_b(\Omega; \mathbb{R}^m)$. If $\mu_k \overset{*}{\rightharpoonup} \mu$ weakly* in $\mathcal{M}_b(\Omega; \mathbb{R}^m)$, then*

$$\int_\Omega H\Big(x, \frac{\mathrm{d}\mu}{\mathrm{d}|\mu|}(x)\Big)\, \mathrm{d}|\mu|(x) \leq \liminf_{k \to +\infty} \int_\Omega H\Big(x, \frac{\mathrm{d}\mu_k}{\mathrm{d}|\mu_k|}(x)\Big)\, \mathrm{d}|\mu_k|(x),$$

*for every lower semicontinuous function $H \colon \Omega \times \mathbb{R}^m \to [0, \infty]$, positively 1-homogeneous and convex in the second variable.*

**Functions of bounded deformation.** Let $\Omega$ be an open set in $\mathbb{R}^n$. For every $u \in L^1(\Omega; \mathbb{R}^n)$, we denote by $\mathrm{E}u$ the $\mathbb{M}^{n \times n}_{\mathrm{sym}}$-valued distribution on $\Omega$, whose components are given by $\mathrm{E}_{ij}u := \frac{1}{2}(D_j u^i + D_i u^j)$. The space $BD(\Omega)$ of *functions of bounded deformation* is the space of all $u \in L^1(\Omega; \mathbb{R}^n)$ such that $\mathrm{E}u \in \mathcal{M}_b(\Omega; \mathbb{M}^{n \times n}_{\mathrm{sym}})$.

A sequence $(u_k)_k$ converges to $u$ weakly* in $BD(\Omega)$ if and only if $u_k \to u$ strongly in $L^1(\Omega; \mathbb{R}^n)$ and $\mathrm{E}u_k \overset{*}{\rightharpoonup} \mathrm{E}u$ weakly* in $\mathcal{M}_b(\Omega; \mathbb{M}^{n \times n}_{\mathrm{sym}})$. We recall that for every $u \in BD(\Omega)$ the measure $\mathrm{E}u$ vanishes on sets of $\mathcal{H}^{n-1}$-measure zero.

The two following embedding theorems hold for the space of functions of bounded deformation. We denote by $1^* := \frac{n}{n-1}$ the Sobolev conjugate of $1$.

**Theorem 2.3.** *The space $BD(\mathbb{R}^n)$ is continuously embedded in $L^{1^*}(\mathbb{R}^n;\mathbb{R}^n)$. More precisely, there exists a constant $C_1 = C_1(n) > 0$ such that for every $u \in BD(\mathbb{R}^n)$ we have*

$$\|u\|_{L^{1^*}(\mathbb{R}^n;\mathbb{R}^n)} \leq C_1 |\mathrm{E}u|(\mathbb{R}^n).$$

*If $\Omega$ is a bounded, open set with Lipschitz boundary, the space $BD(\Omega)$ is continuously embedded in $L^q(\Omega;\mathbb{R}^n)$ for every $1 \leq q \leq 1^*$.*

**Theorem 2.4.** *Let $\Omega$ be a bounded, open set with Lipschitz boundary. Then the space $BD(\Omega)$ is compactly embedded in $L^q(\Omega;\mathbb{R}^n)$ for every $1 \leq q < 1^*$.*

We refer to the book [33] for more details on the general properties of functions of bounded deformation and to [5] for their fine properties.

**Capacity.** For the notion of capacity we refer, e.g., to [17, 19]. We recall here the definition and some properties.

Let $1 \leq q < +\infty$ and let $\Omega$ be a bounded, open subset of $\mathbb{R}^n$. For every subset $B \subset \Omega$, the $q$-*capacity* of $E$ in $\Omega$ is defined by

$$\mathrm{Cap}_q(E,\Omega) := \inf \left\{ \int_\Omega |\nabla v|^q \, \mathrm{d}x \ : \ v \in W_0^{1,q}(\Omega), \ v \geq 1 \text{ a.e. in a neighbourhood of } E \right\}.$$

A set $E \subset \Omega$ has $q$-*capacity zero* if $\mathrm{Cap}_q(E,\Omega) = 0$ (actually, the definition does not depend on the open set $\Omega$ containing $E$). A property is said to hold $\mathrm{Cap}_q$-*quasi everywhere* (abbreviated as $\mathrm{Cap}_q$-q.e.) if it holds for a set of $q$-capacity zero.

If $1 < q \leq n$ and $E$ has $q$-capacity zero, then $\mathcal{H}^s(E) = 0$ for every $s > n - q$.

A function $\alpha \colon \Omega \to \mathbb{R}$ is $\mathrm{Cap}_q$-*quasicontinuous* if for every $\varepsilon > 0$ there exists a set $E_\varepsilon \subset \Omega$ with $\mathrm{Cap}_q(E_\varepsilon) < \varepsilon$ such that the restriction $\alpha|_{\Omega \setminus E_\varepsilon}$ is continuous. Note that if $q > n$, a function $\alpha$ is $\mathrm{Cap}_q$-quasicontinuous if and only if it is continuous.

Every function $\alpha \in W^{1,q}(\Omega)$ admits a $\mathrm{Cap}_q$-*quasicontinuous representative* $\widetilde{\alpha}$, i.e., a $\mathrm{Cap}_q$-quasicontinuous function $\widetilde{\alpha}$ such that $\widetilde{\alpha} = \alpha$ $\mathcal{L}^n$-a.e. in $\Omega$. The $\mathrm{Cap}_q$-quasicontinuous representative is essentially unique, that is, if $\widetilde{\beta}$ is another $\mathrm{Cap}_q$-quasicontinuous representative of $\alpha$, then $\widetilde{\beta} = \widetilde{\alpha}$ $\mathrm{Cap}_q$-q.e. in $\Omega$. If $\alpha_k \to \alpha$ strongly in $W^{1,q}(\Omega)$, then there exists a subsequence $k_j$ such that $\widetilde{\alpha}_{k_j} \to \widetilde{\alpha}$ $\mathrm{Cap}_q$-q.e. in $\Omega$.

## 3. Concentration phenomena

In the whole section we assume that $\Omega$ is a bounded, open set with Lipschitz boundary.

In order to prove the lower semicontinuity result, we shall provide a precise description of the weak* limit of the sequence of measures $\widetilde{\alpha}_k \mathrm{E}u_k$, for $\alpha_k \rightharpoonup \alpha$ weakly in $W^{1,n}(\Omega)$ and $u_k \overset{*}{\rightharpoonup} u$ weakly* in $BD(\Omega)$. We start by showing that, in general, the sequence $\widetilde{\alpha}_k \mathrm{E}u_k$ does not converge to $\widetilde{\alpha} \mathrm{E}u$ weakly* in $\mathcal{M}_b(\Omega;\mathbb{M}^{n \times n}_{\mathrm{sym}})$. Indeed concentration phenomena may occur, as the following example shows.

*Example* 3.1. Let $n = 2$ and let $\Omega = (-1,1)^2$. We construct here an explicit example of a sequence $(\alpha_k)_k$ in $W^{1,2}(\Omega)$ with $0 \leq \alpha_k \leq 1$ and a sequence $(u_k)_k$ in $BD(\Omega)$ such that

$$\alpha_k \rightharpoonup 0 \quad \text{weakly in } W^{1,2}(\Omega), \tag{3.1}$$

$$u_k \overset{*}{\rightharpoonup} 0 \quad \text{weakly* in } BD(\Omega), \tag{3.2}$$

but nonetheless

$$\widetilde{\alpha}_k \mathrm{E}u_k \text{ does not converge to } 0 \text{ weakly* in } \mathcal{M}_b(\Omega;\mathbb{M}^{2 \times 2}_{sym}). \tag{3.3}$$

Let us define the polygon $P_k = A_k \cup B_k \cup C_k \cup D_k$ as in Figure 1. Let $A_k := \left(-\frac{1}{2k}, \frac{1}{2k}\right) \times \left(-\frac{1}{k}, 0\right)$ and $B_k := \left(-\frac{1}{2k}, \frac{1}{2k}\right) \times \left(0, \frac{1}{k}\right)$. Let $C_k$ be the union of the triangle $C_k^+$ with vertices $\left(\frac{1}{2k}, 0\right)$, $\left(\frac{3}{2k}, 0\right)$, $\left(\frac{1}{2k}, \frac{1}{k}\right)$ and of the triangle $C_k^-$ with vertices $\left(-\frac{1}{2k}, 0\right)$, $\left(-\frac{3}{2k}, 0\right)$, $\left(-\frac{1}{2k}, \frac{1}{k}\right)$. Let $D_k$ be the union of the triangle $D_k^+$ with vertices $\left(\frac{1}{2k}, 0\right)$, $\left(\frac{1}{2k}, -\frac{1}{k}\right)$, $\left(\frac{3}{2k}, 0\right)$ and of the triangle $D_k^-$ with vertices $\left(-\frac{1}{2k}, 0\right)$, $\left(-\frac{1}{2k}, -\frac{1}{k}\right)$, $\left(-\frac{3}{2k}, 0\right)$. For $k$ large enough, $P_k$ is contained in $\Omega$.

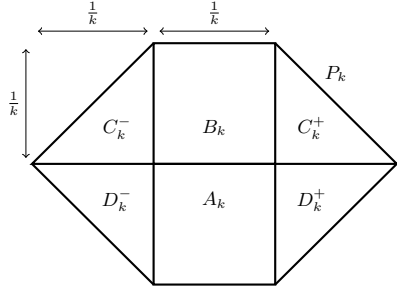

**Figure 1.** Decomposition of the set $P_k$.

We define the piecewise affine functions $\alpha_k \in W^{1,\infty}(\Omega)$ in such a way that $\alpha_k(x) = 1$ for every $x \in \partial A_k \cap \partial B_k = \left(-\frac{1}{2k}, \frac{1}{2k}\right) \times \{0\}$, $\alpha_k(x) = 0$ for every $x \notin P_k$, and $\alpha_k$ is affine on each of the sets which decompose $P_k$. Notice that $0 \leq \alpha_k \leq 1$ and that

$$\nabla \alpha_k(x) = \begin{cases} ke_2 & \text{if } x \in A_k, \\ -ke_2 & \text{if } x \in B_k, \\ -ke_1 - ke_2 & \text{if } x \in C_k^+, \\ -ke_1 + ke_2 & \text{if } x \in D_k^+, \\ ke_1 - ke_2 & \text{if } x \in C_k^-, \\ ke_1 + ke_2 & \text{if } x \in D_k^-, \end{cases}$$

where $\{e_1, e_2\}$ is the standard basis in $\mathbb{R}^2$. In particular $\sup_k \|\nabla \alpha_k\|_{L^2(\Omega;\mathbb{R}^2)} < +\infty$ and $\nabla \alpha_k \to 0$ strongly in $L^1(\Omega; \mathbb{R}^2)$. Finally, we define $u_k \colon \mathbb{R}^2 \to \mathbb{R}^2$ by $u_k := |\nabla \alpha_k| \mathbb{1}_{A_k} e_1 = k \mathbb{1}_{A_k} e_1$, where $\mathbb{1}_{A_k}$ is the indicator function of the set $A_k$.

Since

$$\int_\Omega |\alpha_k(x)|^2 \, dx \leq \mathcal{L}^2(P_k) \to 0,$$

$$\sup_k \int_\Omega |\nabla \alpha_k(x)|^2 \, dx < +\infty,$$

we deduce (3.1). Moreover

$$\int_\Omega |u_k(x)| \, dx = \int_{A_k} |\nabla \alpha_k(x)| \, dx \to 0,$$

$$\sup_k |\mathrm{E} u_k|(\Omega) \leq C \sup_k \left[ k \, \mathcal{H}^1(\partial A_k) \right] < +\infty$$

imply (3.2). In order to prove (3.3), let us fix $\varphi \in \mathcal{C}_0(\Omega; \mathbb{M}^{2 \times 2}_{sym})$. Let us denote the sides of $A_k$ by $L_k^i$, $i = 1, 2, 3, 4$, $L_k^1$ being the upper side and $L_k^3$ being the bottom side. Notice that the measure $\widetilde{\alpha}_k \, \mathrm{E} u_k$ is concentrated on $L_k^1 \cup L_k^2 \cup L_k^4$ and that

$$\int_{L_k^2 \cup L_k^4} \widetilde{\alpha}_k \, \varphi : d\mathrm{E} u_k = k^2 \int_{-\frac{1}{k}}^{0} x_2 \left[ \varphi\left(-\tfrac{1}{2k}, x_2\right) - \varphi\left(\tfrac{1}{2k}, x_2\right) \right] : e_1 \odot e_1 \, dx_2 \to 0,$$

$x_2$ denoting the second coordinate of $x$. Therefore, the only contribution to the limit is given by

$$\int_{L_k^1} \widetilde{\alpha}_k \, \varphi : d\mathrm{E} u_k = -k \int_{L_k^1} \varphi : (e_1 \odot e_2) \, d\mathcal{H}^1 \to -\varphi(0) : (e_1 \odot e_2),$$

i.e., $\widetilde{\alpha}_k \mathrm{E} u_k \overset{*}{\rightharpoonup} -\delta_0 \, e_1 \odot e_2$. This proves the claim. The example can be also modified in order to have $\mathrm{div}\, u_k = 0$. This can be done by suitably extending the vector field $u_k$ in $D_k^+$, $D_k^-$, and $\Omega \setminus P_k$.

In the previous example, the difference between $\widetilde{\alpha} \, \mathrm{E} u = 0$ and the weak* limit of $\widetilde{\alpha}_k \mathrm{E} u_k$ is a measure concentrated on a point. Actually, we will show that for every sequence $(\widetilde{\alpha}_k \mathrm{E} u_k)_k$ the excess measure in the limit may concentrate on at most countably many points. Specifically, we shall prove the following result.

**Theorem 3.2.** *Let $\alpha_k, \alpha \in W^{1,n}(\Omega)$ and $u_k, u \in BD(\Omega)$. Assume that*

$$\|\alpha_k\|_{L^\infty(\Omega)} \leq M,  \tag{3.4}$$

$$\alpha_k \rightharpoonup \alpha \quad \text{weakly in } W^{1,n}(\Omega), \tag{3.5}$$

$$u_k \stackrel{*}{\rightharpoonup} u \quad \text{weakly* in } BD(\Omega). \tag{3.6}$$

*Then, up to a subsequence (which we do not relabel),*

$$\widetilde{\alpha}_k \operatorname{E} u_k \stackrel{*}{\rightharpoonup} \widetilde{\alpha} \operatorname{E} u + \eta \quad \text{weakly* in } \mathcal{M}_b(\Omega; \mathbb{M}^{n\times n}_{\operatorname{sym}}),$$

*where $\eta \in \mathcal{M}_b(\Omega; \mathbb{M}^{n\times n}_{\operatorname{sym}})$ is concentrated on an at most countable set.*

The initial step for the proof of Theorem 3.2 is a careful analysis of the limit behaviour of a sequence $(u_k)_k$ converging weakly* in $BD(\Omega)$. The Embedding Theorems for $BD(\Omega)$ (Theorem 2.3 and Theorem 2.4) do not guarantee that the sequence $(u_k)_k$ converges strongly in $L^{1^*}(\Omega; \mathbb{R}^n)$. Nevertheless, the following concentration compactness argument in the spirit of [21, 22] shows that the lack of compactness of $(u_k)_k$ in $L^{1^*}(\Omega; \mathbb{R}^n)$ is only due to concentration around countably many points. For a proof of the analogous result in the Sobolev case we refer e.g. to [16].

**Theorem 3.3.** *Let $(u_k)_k$ be a sequence in $BD(\Omega)$. Assume that $u_k \stackrel{*}{\rightharpoonup} 0$ weakly* in $BD(\Omega)$ and that*

$$|u_k|^{1^*} \stackrel{*}{\rightharpoonup} \nu \quad \text{weakly* in } \mathcal{M}_b(\Omega) \tag{3.7}$$

*for some non-negative measure $\nu \in \mathcal{M}_b^+(\Omega)$. Then $\nu$ is concentrated on an at most countable set, i.e., there exists a countable set $\{x_j\}_j$ of points of $\Omega$ such that*

$$\nu = \sum_j c_j \delta_{x_j},$$

*with $c_j \in (0, +\infty)$.*

*Proof.* Upon extracting a subsequence (which we do not relabel), we suppose that

$$|\operatorname{E} u_k| \stackrel{*}{\rightharpoonup} \mu \quad \text{weakly* in } \mathcal{M}_b(\Omega) \tag{3.8}$$

for some measure non-negative measure $\mu \in \mathcal{M}_b^+(\Omega)$. Let us define the set

$$D := \{x \in \Omega \, : \, \mu(\{x\}) > 0\}.$$

Note that the set $D$ is at most countable, since $\mu$ is a finite measure. We claim that $\nu$ is concentrated on a subset of $D$.

We first prove that the measure $\nu$ is absolutely continuous with respect to $\mu$. Let us fix a compact set $K \subset \Omega$, and an open set $V \subset \Omega$ such that $K \subset V$. Let us consider a cut-off function $\phi \in \mathcal{C}_c^1(\Omega)$ with $0 \leq \phi \leq 1$, $\phi = 1$ on $K$, $\operatorname{supp}(\phi) \subset V$. The functions $\phi u_k$ have compact support in $\Omega$, they belong to $BD(\mathbb{R}^n)$, and $\operatorname{E}(\phi u_k) = \phi \operatorname{E} u_k + \nabla \phi \odot u_k$. By Theorem 2.3, we infer that

$$\left( \int_{\mathbb{R}^n} |\phi u_k|^{1^*} \mathrm{d}x \right)^{1/1^*} \leq C_1 |\operatorname{E}(\phi u_k)|(\mathbb{R}^n) \leq C_1 \left[ \int_\Omega |\phi| \, \mathrm{d}|\operatorname{E} u_k| + \int_\Omega |\nabla \phi \odot u_k| \, \mathrm{d}x \right].$$

Since $u_k \to 0$ strongly in $L^1(\Omega; \mathbb{R}^n)$ (Theorem 2.4), we have

$$\int_\Omega |\nabla \phi \odot u_k| \, \mathrm{d}x \to 0$$

as $k \to +\infty$. Testing (3.7) and (3.8) with the functions $|\phi|^{1^*}$ and $|\phi|$ respectively, we pass to the limit as $k \to +\infty$ in the inequality above and we get

$$\left( \int_\Omega |\phi|^{1^*} \mathrm{d}\nu \right)^{1/1^*} \leq C_1 \int_\Omega |\phi| \, \mathrm{d}\mu.$$

From the assumptions on $\phi$ we deduce that

$$\left( \nu(K) \right)^{1/1^*} \leq C_1 \mu(V).$$

By the arbitrariness of $K$ and $V$, we have

$$\left( \nu(B) \right)^{1/1^*} \leq C_1 \mu(B) \tag{3.9}$$

for any Borel set $B \subset \Omega$. Therefore we conclude that $\nu$ is absolutely continuous with respect to $\mu$.

By the Radon-Nikodym Theorem
$$\nu = \frac{d\nu}{d\mu}\mu,$$
where $\frac{d\nu}{d\mu}$ is the Radon-Nikodym derivative of $\nu$ with respect to $\mu$ given by
$$\frac{d\nu}{d\mu}(x) = \lim_{r\to 0^+} \frac{\nu(B_r(x))}{\mu(B_r(x))} \quad \text{for } \mu\text{-a.e. } x \in \Omega.$$
By (3.9) and the formula above we infer that
$$\frac{d\nu}{d\mu}(x) \leq \limsup_{r\to 0^+} \left[ C_1^{1^*} \mu(B_r(x))^{1^*-1} \right] = 0 \quad \text{for } \mu\text{-a.e. } x \in \Omega \setminus D,$$
i.e., that $\nu$ is concentrated on a subset of $D$. $\square$

The following lemma will be used in the proof of Theorem 3.2 to characterise the limit of the sequence $(\nabla \alpha_k \odot u_k)_k$.

**Lemma 3.4.** *Let $(g_k)_k$ be a bounded sequence in $L^n(\Omega; \mathbb{R}^n)$ and let $(u_k)_k$ be a sequence in $BD(\Omega)$ such that $u_k \overset{*}{\rightharpoonup} 0$ weakly* in $BD(\Omega)$. Assume that*
$$|g_k \odot u_k| \overset{*}{\rightharpoonup} \nu \quad \text{weakly* in } \mathcal{M}_b(\Omega)$$
*for some non-negative measure $\nu \in \mathcal{M}_b^+(\Omega)$. Then $\nu$ is concentrated on an at most countable set.*

*Proof.* By Theorem 2.3, the sequence $(|u_k|^{1^*})_k$ is bounded in $L^1(\Omega)$. Upon extracting a subsequence (which we do not relabel), we suppose that
$$|g_k|^n \overset{*}{\rightharpoonup} \nu^g, \quad |u_k|^{1^*} \overset{*}{\rightharpoonup} \nu^u \quad \text{weakly* in } \mathcal{M}_b(\Omega).$$
Let us fix a compact set $K \subset \Omega$, and an open set $V \subset \Omega$ such that $K \subset V$. Let $\phi \in \mathcal{C}_c^1(\Omega)$ be such that $0 \leq \phi \leq 1$, $\phi = 1$ on $K$, and $\text{supp}(\phi) \subset V$. By Hölder's Inequality we have
$$\int_\Omega \phi^2 |g_k \odot u_k| \, dx \leq C \int_\Omega |\phi\, g_k||\phi\, u_k| \, dx \leq C \left( \int_\Omega |\phi|^n |g_k|^n \, dx \right)^{1/n} \left( \int_\Omega |\phi|^{1^*} |u_k|^{1^*} \, dx \right)^{1/1^*}.$$
Passing to the limit as $k \to +\infty$ we deduce that
$$\int_\Omega \phi^2 \, d\nu \leq C \left( \int_\Omega |\phi|^n \, d\nu^g \right)^{1/n} \left( \int_\Omega |\phi|^{1^*} \, d\nu^u \right)^{1/1^*}$$
and thus
$$\nu(K) \leq C \big(\nu^g(V)\big)^{1/n} \big(\nu^u(V)\big)^{1/1^*}.$$
By the arbitrariness of $K$ and $V$ we conclude that
$$\nu(B) \leq C \big(\nu^g(B)\big)^{1/n} \big(\nu^u(B)\big)^{1/1^*}$$
for every Borel set $B$, and therefore that $\nu$ is absolutely continuous with respect to $\nu^u$, which by Theorem 3.3 is concentrated on an at most countable set. $\square$

We shall need the following Leibniz rule formula for the product of Sobolev functions and functions of bounded deformation. We include the proof for the convenience of the reader.

**Proposition 3.5.** *Let $\alpha \in W^{1,n}(\Omega) \cap L^\infty(\Omega)$, and $u \in BD(\Omega)$. Then $\alpha u \in BD(\Omega)$ and*
$$\mathrm{E}(\alpha u) = \widetilde{\alpha}\, \mathrm{E}u + \nabla \alpha \odot u \quad \text{in } \mathcal{M}_b(\Omega; \mathbb{M}_{\text{sym}}^{n\times n}). \tag{3.10}$$

*Proof.* The proof is based on an approximation argument. There exists a sequence of smooth functions $\alpha_k \in C^\infty(\overline{\Omega})$ such that $\alpha_k \to \alpha$ strongly in $W^{1,n}(\Omega)$ and $\|\alpha_k\|_{L^\infty(\Omega)} \leq \|\alpha\|_{L^\infty(\Omega)}$. It is immediate to prove via integration by parts that
$$\mathrm{E}(\alpha_k u) = \alpha_k \mathrm{E}u + \nabla \alpha_k \odot u \quad \text{in } \mathcal{M}_b(\Omega; \mathbb{M}_{\text{sym}}^{n\times n}).$$
In particular, the total variations $|\mathrm{E}(\alpha_k u)|$ are bounded, and thus $\mathrm{E}(\alpha_k u) \overset{*}{\rightharpoonup} \mathrm{E}(\alpha u)$. Moreover, $\nabla \alpha_k \odot u \to \nabla \alpha \odot u$ strongly in $L^1(\Omega; \mathbb{M}_{\text{sym}}^{n\times n})$. To conclude the proof of (3.10), we simply remark that $\alpha_k \to \widetilde{\alpha}$ Cap$_n$-q.e. (up to a subsequence) and $\mathrm{E}u$ vanishes on sets of $n$-capacity zero, so that $\alpha_k \mathrm{E}u \overset{*}{\rightharpoonup} \widetilde{\alpha}\, \mathrm{E}u$ in $\mathcal{M}_b(\Omega; \mathbb{M}_{\text{sym}}^{n\times n})$ by the Dominated Convergence Theorem. $\square$

We are now in a position to prove Theorem 3.2.

8  V. CRISMALE AND G. ORLANDO8 V. CRISMALE AND G. ORLANDO

*Proof of Theorem 3.2.* By Proposition 3.5 we have

$$\widetilde{\alpha}_k \operatorname{E} u_k = \operatorname{E}(\alpha_k u_k) - \nabla \alpha_k \odot u_k \,. \tag{3.11}$$

Notice that

$$\operatorname{E}(\alpha_k u_k) \overset{*}{\rightharpoonup} \operatorname{E}(\alpha\, u) \quad \text{weakly* in } \mathcal{M}_b(\Omega; \mathbb{M}^{n\times n}_{\mathrm{sym}})\,. \tag{3.12}$$

Indeed, by Hölder's Inequality

$$|\operatorname{E}(\alpha_k u_k)|(\Omega) \leq \|\alpha_k\|_{L^\infty(\Omega)} |\operatorname{E} u_k|(\Omega) + C \|\nabla \alpha_k\|_{L^n(\Omega; \mathbb{R}^n)} \|u_k\|_{L^{1^*}(\Omega; \mathbb{R}^n)}\,.$$

By (3.4)–(3.6) and by Theorem 2.3 the right-hand side in the inequality above is uniformly bounded. Since $\alpha_k u_k \to \alpha\, u$ strongly in $L^1(\Omega; \mathbb{R}^n)$, we conclude that (3.12) holds.

We now study the weak* limit of $(\nabla \alpha_k \odot u_k)_k$ in $\mathcal{M}_b(\Omega; \mathbb{M}^{n\times n}_{\mathrm{sym}})$. Since $\nabla \alpha_k \rightharpoonup \nabla \alpha$ weakly in $L^n(\Omega; \mathbb{R}^n)$, we get that

$$\nabla \alpha_k \odot u \rightharpoonup \nabla \alpha \odot u \quad \text{weakly in } L^1(\Omega; \mathbb{M}^{n\times n}_{\mathrm{sym}})\,. \tag{3.13}$$

Upon the extraction of a subsequence (that we do not relabel), we can assume that

$$\nabla \alpha_k \odot (u - u_k) \overset{*}{\rightharpoonup} \eta \quad \text{weakly* in } \mathcal{M}_b(\Omega; \mathbb{M}^{n\times n}_{\mathrm{sym}})\,, \tag{3.14}$$

$$|\nabla \alpha_k \odot (u - u_k)| \overset{*}{\rightharpoonup} \nu \quad \text{weakly* in } \mathcal{M}_b(\Omega)\,.$$

By Lemma 3.4 we have that $\nu$, and *a fortiori* $\eta$, is concentrated on an at most countable set.

By (3.13) and (3.14) we get that

$$\nabla \alpha_k \odot u_k \overset{*}{\rightharpoonup} \nabla \alpha \odot u - \eta \tag{3.15}$$

weakly* in $\mathcal{M}_b(\Omega; \mathbb{M}^{n\times n}_{\mathrm{sym}})$.

From (3.11), (3.12), (3.15), and Proposition 3.5 we conclude that

$$\widetilde{\alpha}_k \operatorname{E} u_k \overset{*}{\rightharpoonup} \operatorname{E}(\alpha\, u) - \nabla \alpha \odot u + \eta = \widetilde{\alpha}\operatorname{E} u + \eta$$

weakly* in $\mathcal{M}_b(\Omega; \mathbb{M}^{n\times n}_{\mathrm{sym}})$. $\square$

*Remark* 3.6. Theorem 3.2 does not hold if $\alpha_k \in W^{1,q}(\Omega)$ with $q < n$. In this case, the difference between $\widetilde{\alpha}\operatorname{E} u$ and the weak* limit of $\widetilde{\alpha}_k \operatorname{E} u_k$ may be not singular with respect to measures which vanish on sets with Hausdorff dimension strictly less than $n-1$. We provide an example below.

*Example* 3.7. Let $n = 2$, let $\Omega = (-2, 2)^2$, and let $1 < q < 2$. We provide here an example of a sequence $(\beta_k)_k$ in $W^{1,q}(\Omega)$ with $0 \leq \beta_k \leq 1$ and a sequence $(u_k)_k$ in $BD(\Omega)$ such that $\beta_k \rightharpoonup 0$ weakly in $W^{1,q}(\Omega)$, $u_k \overset{*}{\rightharpoonup} 0$ weakly* in $BD(\Omega)$, and the weak* limit of $\widetilde{\beta}_k \operatorname{E} u_k$ is concentrated on a set of Hausdorff dimension 1.

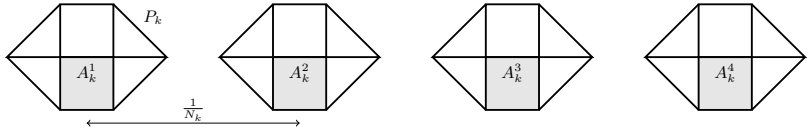

**Figure 2.** The function $\beta_k$ is supported on the union of the $N_k$ equispaced copies of $P_k$, while the function $u_k$ is supported on the grey region given by $\bigcup_{j=1}^{N_k} A_k^j$.

Let $\alpha_k \in W^{1,\infty}(\Omega)$ be the piecewise affine functions supported on the polygons $P_k$ and let $A_k$ be the cubes exhibited in Example 3.1. Let $N_k$ be the integer part of $k^{2-q}$ and let $x_k^j = \left(\frac{j-1}{N_k}, 0\right)$. We define $\beta_k(x) := \sum_{j=1}^{N_k} \alpha_k(x - x_k^j)$ and $u_k := \sum_{j=1}^{N_k} k^{q-1} e_1 \mathbb{1}_{A_k^j}$, where $A_k^j = A_k + x_k^j$. (See Figure 2.) Notice that $\beta_k \to 0$ strongly in $L^q(\Omega)$, $u_k \to 0$ strongly in $L^1(\Omega; \mathbb{R}^2)$, and

$$\int_\Omega |\nabla \beta_k|^q \, \mathrm{d}x = \sum_{j=1}^{N_k} \int_{P_k} |\nabla \alpha_k|^q \, \mathrm{d}x \sim N_k \frac{1}{k^2} k^q \sim 1\,,$$

$$|\operatorname{E} u_k|(\Omega) \leq C \sum_{j=1}^{N_k} k^{q-1} \mathcal{H}^1(\partial A_k) \sim N_k\, k^{q-1} \frac{1}{k} \sim 1\,,$$

as $k \to +\infty$. Thus $\beta_k \rightharpoonup 0$ weakly in $W^{1,q}(\Omega)$ and $u_k \overset{*}{\rightharpoonup} 0$ weakly* in $BD(\Omega)$. With computations similar to those contained in Example 3.1, it is easy to show that only the restriction of $\widetilde{\beta}_k \operatorname{E} u_k$ to the upper sides of the squares $A_k^j$ gives a contribution to the limit. Hence for every $\varphi \in \mathcal{C}_0(\Omega; \mathbb{M}_{sym}^{2\times 2})$ we have

$$\int_\Omega \widetilde{\beta}_k \, \varphi : \mathrm{dE}u_k = -\sum_{j=1}^{N_k} k^{q-1} \int_{-\frac{1}{2k}}^{\frac{1}{2k}} \varphi\big(x_1 + \tfrac{j-1}{N_k}, 0\big) : e_1 \odot e_2 \, \mathrm{d}x_1 + o(1) \to - \int_{[0,1]\times\{0\}} \varphi : e_1 \odot e_2 \, \mathrm{d}\mathcal{H}^1,$$

i.e., $\widetilde{\beta}_k \, \mathrm{E}u_k \overset{*}{\rightharpoonup} -e_1 \odot e_2 \mathcal{H}^1 \llcorner \big([0,1]\times\{0\}\big)$.

## 4. Proof of Theorem 1.1

Upon the extraction of a subsequence (that we do not relabel), we assume that the liminf in (1.5) is actually a limit.

We shall prove the theorem supposing that $V$ is a Lipschitz function. Indeed, if this is not the case, we can always find an increasing family of Lipschitz functions $V_h \colon \mathbb{R} \to [0, +\infty)$ such that $V = \sup_h V_h$. Then, assuming that (1.5) holds for each $V_h$, we have

$$\int_\Omega V_h(\widetilde{\alpha}(x)) \, H\Big(x, \frac{\mathrm{d}p}{\mathrm{d}|p|}(x)\Big) \, \mathrm{d}|p|(x) \leq \liminf_{k\to+\infty} \int_\Omega V_h(\widetilde{\alpha}_k(x)) \, H\Big(x, \frac{\mathrm{d}p_k}{\mathrm{d}|p_k|}(x)\Big) \, \mathrm{d}|p_k|(x)$$

$$\leq \liminf_{k\to+\infty} \int_\Omega V(\widetilde{\alpha}_k(x)) \, H\Big(x, \frac{\mathrm{d}p_k}{\mathrm{d}|p_k|}(x)\Big) \, \mathrm{d}|p_k|(x),$$

and by the Monotone Convergence Theorem we deduce (1.5).

Let us define the non-negative functions $\beta_k := V(\alpha_k)$ and $\beta := V(\alpha)$. Since $V$ is Lipschitz and $\Omega$ is bounded, the chain rule for Sobolev functions implies that $\beta_k, \beta \in W^{1,n}(\Omega)$. Moreover, it is immediate to see that $\beta_k \rightharpoonup \beta$ weakly in $W^{1,n}(\Omega)$, i.e., the sequence $(\beta_k)_k$ satisfies the same assumptions on the sequence $(\alpha_k)_k$. Moreover, $\beta_k \geq 0$ a.e. in $\Omega$.

Let us prove the theorem under the additional assumption that $\|\beta_k\|_{L^\infty(\Omega)} \leq M$. Notice that $\beta_k \to \beta$ strongly in $L^2(\Omega)$. Together with (1.4), this implies that $\beta_k e_k \rightharpoonup \beta e$ weakly in $L^1(\Omega; \mathbb{M}_{sym}^{n\times n})$. Hence, by (1.2) and Theorem 3.2, we have (up to a subsequence)

$$\widetilde{\beta}_k \, p_k = \widetilde{\beta}_k \, \mathrm{E}u_k - \beta_k \, e_k \overset{*}{\rightharpoonup} \widetilde{\beta} \, \mathrm{E}u - \beta \, e + \eta = \widetilde{\beta} \, p + \eta \quad \text{weakly* in } \mathcal{M}_b(\Omega; \mathbb{M}_{sym}^{n\times n}),$$

where the measure $\eta \in \mathcal{M}_b(\Omega; \mathbb{M}_{sym}^{n\times n})$ is concentrated on an at most countable set. Since $|p|$ is concentrated on sets of dimension at most $n-1$, the measures $|\widetilde{\beta} p|$ and $|\eta|$ are mutually singular. By Remark 2.1, by the 1-homogeneity of $H$, and by Reshetnyak's Lower Semicontinuity Theorem we infer that

$$\int_\Omega \widetilde{\beta} \, H\Big(x, \frac{\mathrm{d}p}{\mathrm{d}|p|}\Big) \mathrm{d}|p| \leq \int_\Omega \widetilde{\beta} \, H\Big(x, \frac{\mathrm{d}p}{\mathrm{d}|p|}\Big) \mathrm{d}|p| + \int_\Omega H\Big(x, \frac{\mathrm{d}\eta}{\mathrm{d}|\eta|}\Big) \mathrm{d}|\eta|$$

$$= \int_\Omega H\Big(x, \frac{\mathrm{d}(\widetilde{\beta} p + \eta)}{\mathrm{d}|\widetilde{\beta} p + \eta|}\Big) \mathrm{d}|\widetilde{\beta} p + \eta| \leq \liminf_{k\to+\infty} \int_\Omega \widetilde{\beta}_k \, H\Big(x, \frac{\mathrm{d}p_k}{\mathrm{d}|p_k|}\Big) \mathrm{d}|p_k|.$$

To remove the assumption that the sequence $(\beta_k)_k$ is bounded in $L^\infty(\Omega)$ we use a truncation argument. For every $M > 0$ we define the functions $\beta_k^M := \beta_k \wedge M$ and $\beta^M := \beta \wedge M$. Since $\beta_k^M \rightharpoonup \beta^M$ weakly in $W^{1,n}(\Omega)$, by the previous step we have

$$\int_\Omega \widetilde{\beta}^M \, H\Big(x, \frac{\mathrm{d}p}{\mathrm{d}|p|}\Big) \mathrm{d}|p| \leq \liminf_{k\to+\infty} \int_\Omega \widetilde{\beta}_k \, H\Big(x, \frac{\mathrm{d}p_k}{\mathrm{d}|p_k|}\Big) \mathrm{d}|p_k|.$$

We conclude applying the Monotone Convergence Theorem as $M \to +\infty$.

## 5. Application to a model for linearised elasto-plasticity coupled with damage

In this section we apply Theorem 1.1 to show the esistence of *energetic solutions* (cf. [25]) for a model which couples small-strain plasticity and damage in $W^{1,n}(\Omega)$ (recall $\Omega \subset \mathbb{R}^n$). The mechanical framework for this coupling has been proposed and analysed in [2, 3] (for further contribution in this direction see, e.g., [31, 32, 30, 1]). The existence of quasistatic evolutions has been proven in [8, 11] via the *energetic approach* and via *vanishing viscosity*, respectively (see e.g. [26] for details and comparison for the two approaches). The notion of quasistatic evolution we give below is similar to the one in [8]. In that paper, the damage variable belongs to $W^{1,q}(\Omega)$, with $q > n$, and in particular it is continuous.

We assume that $\Omega$ is an open bounded set with Lipschitz boundary partitioned as $\partial\Omega = \partial_D\Omega \cup \partial_N\Omega \cup N$, with $\partial_D\Omega$ and $\partial_N\Omega$ relatively open, $\partial_D\Omega \cap \partial_N\Omega = \emptyset$, $\mathcal{H}^{n-1}(N) = 0$, and $\partial_D\Omega \neq \emptyset$.

Moreover, we assume that the common boundary between $\partial_D\Omega$ and $\partial_N\Omega$ is smooth enough, more precisely that [11, (2.2)] holds; this is only needed to ensure a suitable integration by parts formula in the stress-strain duality. Let $[0,T]$ be the time interval where we study the evolution, and $u_D \in AC([0,T]; H^1(\mathbb{R}^n; \mathbb{R}^n))$ be a prescribed Dirichlet datum for the displacement on $\partial_D\Omega$. For simplicity of notation, both the surface forces on $\partial_N\Omega$ and the volume forces are null.

Let us now briefly recall the energetic and dissipative terms involved in the definition of energetic solutions for the present model, referring to [8] for more details.

The *elastic energy* is defined on $L^1(\Omega; [0,1]) \times L^2(\Omega; \mathbb{M}^{n \times n}_{\mathrm{sym}})$ by

$$\mathcal{Q}(\alpha, e) := \frac{1}{2} \int_\Omega \mathbb{C}(\alpha(x))\, e(x) : e(x) \, \mathrm{d}x\,.$$

The *elasticity tensor* $\mathbb{C}(\alpha)$ is a symmetric fourth order tensor for any $\alpha$, Lipschitz and non-decreasing in $\alpha$, equicontinuous and equicoercive with respect to $\alpha$, and it induces a linear map on $\mathbb{M}^{n \times n}_{\mathrm{sym}}$ that preserves the space of symmetric deviatoric matrices $\mathbb{M}^{n \times n}_D$, as well as its orthogonal space $\mathbb{R}I$.

The *plastic potential* is defined on $W^{1,n}(\Omega; [0,1]) \times \mathcal{M}_b(\Omega \cup \partial_D\Omega; \mathbb{M}^{n \times n}_D)$ by

$$\mathcal{H}(\alpha, p) := \int_{\Omega \cup \partial_D\Omega} V(\widetilde\alpha(x))\, H\Big(\frac{\mathrm{d}p}{\mathrm{d}|p|}(x)\Big)\, \mathrm{d}|p|(x)\,.$$

We assume that the function $V \colon [0,1] \to [c_1, \infty)$ is Lipschitz and non-decreasing, and that $c_1 > 0$; $H \colon \mathbb{M}^{n \times n}_D \to [0, \infty)$ is positively 1-homogeneous and convex, with $r|\xi| \leq H(\xi) \leq R|\xi|$, for some $r > 0$. Notice that every $\alpha \in W^{1,n}(\Omega)$ is well defined in $\overline\Omega$ up to a set of $n$-capacity zero, by considering any $W^{1,n}$ extension of $\alpha$ to a larger set $\Omega'$. We remark that the hypoteses on $H$ in [8] are slightly more general (see [8, (2.11)]), here we are in the setting of [8, Remark 2.1].

The *plastic dissipation* in a time interval $[s,t]$ is defined for any $\alpha \colon [s,t] \to W^{1,n}(\Omega; [0,1])$ and any $p \colon [s,t] \to \mathcal{M}_b(\Omega \cup \partial_D\Omega; \mathbb{M}^{n \times n}_D)$ by

$$\mathcal{V}_\mathcal{H}(\alpha, p; s, t) := \sup\Big\{ \sum_{j=1}^N \mathcal{H}(\alpha(t_j), p(t_j) - p(t_{j-1})) \ : \ s = t_0 < t_1 < \cdots < t_N = t,\, N \in \mathbb{N} \Big\}. \quad (5.1)$$

Moreover, we consider a non-negative, continuous, and non-increasing function $d$ and we introduce the functional $D \colon L^1(\Omega; [0,1]) \to [0, \infty)$ defined by $D(\alpha) := \int_\Omega d(\alpha(x))\, \mathrm{d}x$. This term accounts for the energy dissipated during the damage process. For a given $w \in H^1(\mathbb{R}^n; \mathbb{R}^n)$, the set of *admissible plasticity triples for $w$* is

$$A(w) := \{(u, e, p) \in BD(\Omega) \times L^2(\Omega; \mathbb{M}^{n \times n}_{\mathrm{sym}}) \times \mathcal{M}_b(\Omega \cup \partial_D\Omega; \mathbb{M}^{n \times n}_D) :$$
$$\mathrm{E}u = e + p \text{ in } \Omega,\, p = (w - u) \odot \nu\, \mathcal{H}^{n-1} \text{ on } \partial_D\Omega\}.$$

We are now ready to give the definition of *energetic solutions* (or *globally stable quasistatic evolutions*) driven by the boundary datum $u_D$.

**Definition 5.1.** An *energetic solution* is a function $t \mapsto (\alpha(t), u(t), e(t), p(t))$ from $[0,T]$ into $W^{1,n}(\Omega; [0,1]) \times BD(\Omega) \times L^2(\Omega; \mathbb{M}^{n \times n}_{\mathrm{sym}}) \times \mathcal{M}_b(\Omega \cup \partial_D\Omega; \mathbb{M}^{n \times n}_D)$ such that $(u(t), e(t), p(t)) \in A(u_D(t))$ for every $t \in [0, T]$ and the following conditions are satisfied:

(QS0) *irreversibility*: $t \in [0, T] \mapsto \alpha(t, x)$ is non-increasing for every $x \in \Omega$;
(QS1) *global stability*: for any $t \in [0, T]$ and any $\hat\alpha \leq \alpha(t)$, $(\hat u, \hat e, \hat p) \in A(u_D(t))$

$$\mathcal{Q}(\alpha(t), e(t)) + D(\alpha(t)) + \int_\Omega |\nabla \alpha(t; x)|^n \, \mathrm{d}x \leq \mathcal{Q}(\hat\alpha, \hat e) + D(\hat\alpha) + \int_\Omega |\nabla \hat\alpha(x)|^n \, \mathrm{d}x + \mathcal{H}(\hat\alpha, \hat p - p(t));$$

(QS2) *energy balance*: for any $t \in [0, T]$

$$\mathcal{Q}(\alpha(t), e(t)) + D(\alpha(t)) + \int_\Omega |\nabla \alpha(t; x)|^n \, \mathrm{d}x + \mathcal{V}_\mathcal{H}(\alpha, p; 0, t) = \mathcal{Q}(\alpha(0), e(0)) + D(\alpha(0))$$
$$+ \int_\Omega |\nabla \alpha(0; x)|^n \, \mathrm{d}x + \int_0^t \int_\Omega \mathbb{C}(\alpha(s; x)) e(s; x) : \mathrm{E}\dot u_D(s; x)\, \mathrm{d}x\, \mathrm{d}s\,.$$

Thanks to Theorem 1.1, we can prove the following existence result.

**Theorem 5.2.** *Let $\alpha_0 \in W^{1,n}(\Omega; [0,1])$ and $(u_0, e_0, p_0) \in A(u_D(0))$ satisfying the global stability condition (QS1) at the initial time. Then there exists an energetic solution such that $\alpha(0) = \alpha_0$, $u(0) = u_0$, $e(0) = e_0$, $p(0) = p_0$.*

*Remark* 5.3. The definition of evolutions above differs from the one in [8] not only for the damage regularisation. Indeed, in [8] there is a parameter $\lambda \in [0,1]$ that accounts for the interplay between damage growth and cumulation of plastic strain, thus for a fatigue phenomenon. We stated, for simplicity of notation, Definition 5.1 only for the case $\lambda = 0$; one can follow the argument in [8] to prove existence of energetic solutions corresponding to any $\lambda$. Moreover, in [8] the Dirichlet boundary was the whole $\partial\Omega$. As observed in [8], it is a minor point to consider also external volume and surface forces.

*Proof of Theorem 5.2.* We can closely follow the proof of [8, Theorem 4.3], basing on a time incremental approach, which is by now well consolidated. This consists in solving incremental minimisation problems, to obtain discrete-time evolutions that satisfy a discrete global stability and a discrete energy inequality, then passing these conditions to the limit as the time discretisation step tends to 0, to eventually get the energy balance by the global stability. We need lower semicontinuity of $\mathcal{H}$ in order to prove the existence of minimisers for the incremental minimisation problems, and to show the lower semicontinuity of the plastic dissipation, which is a supremum of suitable plastic potentials, as the time discretisation step tends to 0.

The lower semicontinuity of $\mathcal{H}$ is deduced by Theorem 1.1 in the following way. Let $U \subset \mathbb{R}^n$ be a bounded, open, Lipschitz set such that $U \cap \partial\Omega = \partial_D\Omega$. Let $\Omega^* := \Omega \cup U$ and let us define for any $w \in H^1(\mathbb{R}^n; \mathbb{R}^n)$ and for any $(u, e, p) \in A(w)$

$$u^* := \begin{cases} u & \text{in } \Omega, \\ w & \text{in } \Omega^* \setminus \Omega, \end{cases} \qquad e^* := \begin{cases} e & \text{in } \Omega, \\ Ew & \text{in } \Omega^* \setminus \Omega, \end{cases} \qquad p^* := \begin{cases} p & \text{in } \overline{\Omega}, \\ 0 & \text{in } \Omega^* \setminus \overline{\Omega}. \end{cases} \qquad (5.2)$$

We also consider a continuous extension operator from $W^{1,n}(\Omega)$ to $W^{1,n}(\Omega^*)$ and we associate to any $\alpha \in W^{1,n}(\Omega)$ its extension $\alpha^* \in W^{1,n}(\Omega^*)$. Then

$$\mathcal{H}(\alpha, p) = \int_{\Omega^*} V(\widetilde{\alpha}^*(x)) \, H\Big(\frac{\mathrm{d}p^*}{\mathrm{d}|p^*|}(x)\Big) \, \mathrm{d}|p^*|(x) \,.$$

If $\alpha_k \rightharpoonup \alpha$ weakly in $W^{1,n}(\Omega)$, $w_k \rightharpoonup w$ weakly in $H^1(\mathbb{R}^n; \mathbb{R}^n)$, $(u_k, e_k, p_k) \in A(w_k)$, $u_k \overset{*}{\rightharpoonup} u$ weakly* in $BD(\Omega)$, and $e_k \rightharpoonup e$ weakly in $L^2(\Omega; \mathbb{M}^{n \times n}_{\text{sym}})$, by Theorem 1.1 we get

$$\mathcal{H}(\alpha, p) \leq \liminf_{k \to \infty} \mathcal{H}(\alpha_k, p_k) \,.$$

Indeed, by (5.2), we have the convergence $p_k^* \overset{*}{\rightharpoonup} p^*$ weakly* in $\mathcal{M}_b(\Omega^*; \mathbb{M}^{n \times n}_D)$ for the extensions.

With the lower semicontinuity property above at hand, one follows the proof of [8, Theorem 4.3] and concludes Theorem 5.2. In particular, the lower semicontinuity of the plastic dissipation $V_{\mathcal{H}}$ as the time discretisation step tends to 0 follows by the definition (5.1) as supremum of a family of plastic potentials. We remark that, as in [8], no continuity in time of $\alpha$ is required: indeed, the monotonicity in time of $\alpha$ guarantees that the supremum in (5.1) is actually a limit as the maximum step of the partition tends to 0 (cf. [8, Lemma A.1]). This is crucial to deduce the energy balance from the global stability. □

**Acknowledgements.** Vito Crismale has been supported by a public grant as part of the *Investissement d'avenir* project, reference ANR-11-LABX-0056-LMH, LabEx LMH, and is a member of the Gruppo Nazionale per l'Analisi Matematica, la Probabilità e le loro Applicazioni (GNAMPA) of the Istituto Nazionale di Alta Matematica (INdAM). Gianluca Orlando's research has been supported through the TUM University Foundation Fellowship (TUFF).

## References

[1] R. Alessi, J.-J. Marigo, C. Maurini, S. Vidoli. Coupling damage and plasticity for a phase-field regularisation of brittle, cohesive and ductile fracture: One-dimensional examples. *International Journal of Mechanical Sciences*. Published online, doi: 10.1016/j.ijmecsci.2017.05.047.

[2] R. Alessi, J.-J. Marigo, S. Vidoli. Gradient damage models coupled with plasticity and nucleation of cohesive cracks. *Arch. Ration. Mech. Anal.* **214** (2014), 575–615.

[3] R. Alessi, J.-J. Marigo, S. Vidoli. Gradient damage models coupled with plasticity: Variational formulation and main properties. *Mech. Materials* **80** (2015), 351–367.

[4] M. Ambati, R. Kruse, L. De Lorenzis. A phase-field model for ductile fracture at finite strains and its experimental verification. *Comput. Mech.* **57** (2016), 149–167.

[5] L. Ambrosio, A. Coscia, G. Dal Maso. Fine properties of functions with bounded deformation. *Arch. Rational Mech. Anal.* **139** (1997), 201–238.

[6] L. Ambrosio, N. Fusco, D. Pallara. *Functions of bounded variation and free discontinuity problems.* Oxford Mathematical Monographs, The Clarendon Press, Oxford University Press, New York 2000.

[7] J.-F. BABADJIAN, G.A. FRANCFORT, M.G. MORA. Quasistatic evolution in non-associative plasticity the cap model. *SIAM J. Math. Anal.* **44** (2012), 245–292.

[8] V. CRISMALE. Globally stable quasistatic evolution for a coupled elastoplastic-damage model. *ESAIM Control Optim. Calc. Var.* **22** (2016), 883–912.

[9] V. CRISMALE. Globally stable quasistatic evolution for strain gradient plasticity coupled with damage. *Ann. Mat. Pura Appl.* **196** (2017), 641–685.

[10] V. CRISMALE. Some results on quasistatic evolution problems for unidirectional processes. *Ph.D. Thesis* (2016).

[11] V. CRISMALE, G. LAZZARONI. Viscous approximation of quasistatic evolutions for a coupled elastoplastic-damage model. *Calc. Var. Partial Differential Equations* **55**:17 (2016).

[12] G. DAL MASO, A. DESIMONE, M.G. MORA. Quasistatic evolution problems for linearly elastic - perfectly plastic materials. *Arch. Ration. Mech. Anal.* **180** (2006), 237–291.

[13] G. DAL MASO, A. DESIMONE, F. SOLOMBRINO. Quasistatic evolution for Cam-Clay plasticity: a weak formulation via viscoplastic regularization and time rescaling. *Calc. Var. Partial Differential Equations* **40** (2011), 125–181.

[14] G. DAL MASO, G. ORLANDO, R. TOADER. Fracture models for elasto-plastic materials as limits of gradient damage models coupled with plasticity: the antiplane case. *Calc. Var. Partial Differential Equations*, **55**:45 (2016).

[15] G. DAL MASO, G. ORLANDO, R. TOADER. Lower semicontinuity of a class of integral functionals on the space of functions of bounded deformation. *Adv. Cal. Var.* Published online, doi: 10.1515/acv-2015-0036.

[16] L.C. EVANS. *Weak Convergence Methods for Nonlinear Partial Differential Equations.* CBMS 74, American Mathematical Society, 1990.

[17] L.C. EVANS, R.F. GARIEPY *Measure theory and fine properties of functions.* Studies in Advanced Mathematics. CRC Press, Boca Raton, FL, 1992.

[18] G.A. FRANCFORT, M.G. MORA. Quasistatic evolution in non-associative plasticity revisited. *Preprint 2017.*

[19] J. HEINONEN, T. KILPELÄINEN, O. MARTIO. *Nonlinear potential theory of degenerate elliptic equations.* Dover Publications, Inc., Mineola, NY, 2006. Unabridged republication of the 1993 original.

[20] D. KNEES, R. ROSSI, C. ZANINI. A vanishing viscosity approach to a rate-independent damage model. *Math. Models Methods Appl. Sci.*, **23** (2013), 565-616.

[21] P.-L. LIONS. The concentration-compactness principle in the calculus of variations. The limit case, Part I. *Rev. Mat. Iberoamericana*, **1** (1985), No. 1, 145–201.

[22] P.-L. LIONS. The concentration-compactness principle in the calculus of variations. The limit case, Part II. *Rev. Mat. Iberoamericana*, **1** (1985), No. 2, 45–121.

[23] C. MIEHE, F. ALDAKHEEL. A. RAINA. Phase field modeling of ductile fracture at finite strains: A variational gradient-extended plasticity-damage theory. *Internat. J. Plast.* **84** (2016), 1–32.

[24] C. MIEHE, M. HOFACKER, L. SCHÄNZEL, F. ALDAKHEEL. Phase Field Modeling of Fracture in Multi-Physics Problems. Part II. Coupled brittle-to-ductile failure criteria and crack propagation in thermo-elastic-plastic solids. *Comput. Methods Appl. Mech. Engrg.* **294** (2015), 486–522.

[25] A. MIELKE. Evolution of rate-independent systems, in Evolutionary equations. *Vol. II, Handb. Differ. Equ.*, Elsevier/North-Holland, Amsterdam, (2005), 461–559

[26] A. MIELKE, T. ROUBÍČEK. Rate Independent Systems: Theory and Application. Springer, New York, 2015.

[27] K. PHAM, J.-J. MARIGO. Approche variationnelle de lendommagement: I. Les concepts fondamentaux. *Comptes Rendus Mécanique.* **338** (2010), 191-198.

[28] K. PHAM, J.-J. MARIGO. Approche variationnelle de lendommagement: II. Les modèles à gradient. *Comptes Rendus Mécanique.* **338** (2010), 199-206.

[29] Y.G. RESHETNYAK. Weak convergence of completely additive vector functions on a set. *Siberian Math. J.*, **9** (1968), 1039–1045.

[30] R. ROSSI. Existence results for a coupled viscoplastic-damage model in thermoviscoelasticity. *Discrete Contin. Dyn. Syst. Ser. S*, **10** (2017), 1413-1466.

[31] T. ROUBÍČEK, J. VALDMAN. Perfect plasticity with damage and healing at small strains, its modelling, analysis, and computer implementation. *SIAM J. Appl. Math.*, **75** (2016),314-340.

[32] T. ROUBÍČEK, J. VALDMAN. Stress-driven solution to rate-independent elasto-plasticity with damage at small strains and its computer implementation. *Mathematics and Mechanics of Solids*, **22** (2017), 1267–1287.

[33] R. TEMAM. *Problèmes mathématiques en plasticité.* Méthodes Mathématiques de l'Informatique [Mathematical Methods of Information Science], 12. Gauthier-Villars, Montrouge, 1983.

CMAP, ÉCOLE POLYTECHNIQUE, 91128 PALAISEAU CEDEX, FRANCE
*E-mail address*, Vito Crismale: `vito.crismale@polytechnique.edu`

TUM, BOLTZMANNSTRASSE 3, 85747 GARCHING BEI MÜNCHEN, GERMANY
*E-mail address*, Gianluca Orlando: `orlando@ma.tum.de`